\documentclass[10pt,reqno]{amsart}
\usepackage{graphicx}
\theoremstyle{plain}
\newtheorem{thm}{Theorem}[section]
\newtheorem{conj}[thm]{Conjecture}

\title[]{On the leading coefficient of polynomials orthogonal over domains with corners}

\author{Erwin Mi\~{n}a-D\'{\i}az}
\address{University of Mississippi, Department of Mathematics, Hume Hall 305, P. O. Box 1848,  University, MS  38677-1848, USA}
 \email{minadiaz@olemiss.edu}

\date{}
\begin{document}
\begin{abstract}
Let $G$ be the interior domain of a piecewise analytic Jordan curve without cusps. Let $\{p_n\}_{n=0}^\infty$ be the  sequence of polynomials that are orthonormal over $G$ with respect to the area measure, with each $p_n$ having leading coefficient $\lambda_n>0$. It has been proven in \cite{NS} that the asymptotic behavior of $\lambda_n$ as $n\to\infty$ is given by
\[
\frac{n+1}{\pi}\frac{\gamma^{2n+2}}{ \lambda_n^{2}}=1-\alpha_n,
\]
where $\alpha_n=O(1/n)$ as $n\to\infty$ and $\gamma$ is the reciprocal of the logarithmic capacity of the boundary $\partial G$. In this paper, we prove that the $O(1/n)$ estimate for the error term $\alpha_n$ is, in general,  best possible, by exhibiting an example for which  
\[
\liminf_{n\to\infty}\,n\alpha_n>0.
\]
The proof makes use of the Faber polynomials, about which a conjecture is formulated.
\end{abstract}


\maketitle

\section{Introduction}
Let $L$ be a Jordan curve in the complex plane $\mathbb{C}$. The bounded and unbounded components of  $\overline{\mathbb{C}}\setminus L$ will be denoted by $G$ and $\Omega$, respectively. Let $\{p_n\}_{n=0}^\infty$ be the sequence of orthonormal polynomials with respect to the area measure over $G$. That is, each $p_n(z)=\lambda_nz^n+\cdots$ is a polynomial of degree $n$, having positive leading coefficient $\lambda_n$, and for every pair of non-negative integers $m,n,$ we have
\[
\int_G p_n(z)\overline{p_m(z)}dxdy=\delta_{n,m}.
\]

The asymptotic behavior as $n\to\infty$ of these polynomials has been thoroughly investigated when $L$ is an analytic Jordan curve in \cite{Carleman,DMN,DM1,DM2,M1}, while for $L$ having some degree of smoothness, strong asymptotics for $p_n$, outside and on the curve itself, were obtained by Suetin \cite{Suetin}. 

For $L$ a piecewise analytic Jordan curve, investigations on the $n$th-root asymptotics and zero distribution of the polynomials $p_n$ have been carried out in \cite{Levin,May,MSS}. More recently, finer results have been obtained by N. Stylianopoulos in \cite{NS} with the use of  some tools from quasiconformal mapping theory. For instance, it is proven in \cite[Thm. 1.1]{NS} that if $L$ is a piecewise analytic curve without cusps, then the leading coefficients $\lambda_n$ satisfy the asymptotic formula  
\[
\frac{n+1}{\pi}\frac{\gamma^{2n+2}}{ \lambda_n^{2}}=1-\alpha_n
\]
where $\alpha_n=O(1/n)$ as $n\to\infty$, and $\gamma$ is the reciprocal of the logarithmic capacity of $L$.  This quantity $\gamma$ can be introduced in this context via the  conformal map 
\[
\phi:\Omega\to \{w:|w|>1\}
\]
of $\Omega$ onto the exterior of the unit circle, uniquely determined by the conditions $\phi(\infty)=\infty$ and $\phi'(\infty):=\lim_{z\to\infty}\phi(z)/z>0$. This limit is precisely the value of $\gamma$.   

 We notice that \cite{NS} also establishes a strong asymptotic formula for $p_n$ on the exterior of $L$, and several other important estimates and relations that we do not mention here. 
  
In this paper, we prove that the $O(1/n)$ estimate for the error term $\alpha_n$ is, in general,  best possible, by exhibiting an example of a curve $L$ for which  
\[
\liminf_{n\to\infty}\,n\alpha_n>0.
\]

 For each integer $n\geq 0$, the Faber polynomial $F_n$ associated to $L$ \cite{Suetin2} is defined to be the polynomial part of the Laurent expansion  at $\infty$ of $\phi^n$. The polynomials $F_n$ and the functions
\[
E_n(z):=\phi^n(z)-F_n(z),\qquad z\in \Omega, \quad  n\geq 0,
\]
play an important role in the estimation of $\alpha_n$, since, as proven in \cite[Lem. 2.4]{NS},  
\begin{equation}\label{equation3}
\alpha_n=\frac{(n+1)}{\pi}\left\|\frac{F'_{n+1}}{n+1}-\frac{\gamma^{n+1}}{\lambda_n}p_{n}\right\|^2_{L^2(G)}+\frac{\|E'_{n+1}\|^2_{L^2(\Omega)}}{\pi(n+1)}.
\end{equation}

The proof  that $\alpha_n=O(1/n)$ is then accomplished  by showing that \cite[Thm. 2.1]{NS} the first summand  in the right-hand side of (\ref{equation3}) is a big O of the second one, while for the second summand  \cite[Thm. 2.4]{NS} we have
\begin{equation}\label{equation121}
\|E'_{n+1}\|^2_{L^2(\Omega)}=O(1),\qquad (n\to\infty).
\end{equation} 
Summarizing, there exists some constant $C$ independent of $n$ for which 
\begin{equation}\label{equation120}
\frac{\|E'_{n+1}\|^2_{L^2(\Omega)}}{\pi(n+1)}\leq \alpha_n\leq C\frac{\|E'_{n+1}\|^2_{L^2(\Omega)}}{\pi(n+1)}, \quad n\geq 0.
\end{equation}

Consider the circles $C_1=\left\{z:|z-i|=\sqrt{2}\right\}$ and $C_2=\left\{z:|z+i|=\sqrt{2}\right\}$, which intersect at the points $\pm 1$. Let us take $L$ to be the curve consisting of the two arcs of these circles that lie exterior to each other, that is,
\begin{equation}\label{equation1}
L:=\{z\in C_1: \Im(z)\geq 0\}\cup \{z\in C_2: \Im(z)\leq 0\}.
\end{equation}
This is a piecewise analytic curve with corners at $\pm 1$ and exterior angles $\pi/2$. We shall prove the following result.
\begin{thm} \label{thm1}For the curve $L$ defined by (\ref{equation1}), we have
\begin{equation}\label{equation2}
\liminf_{n\to\infty}n\alpha_n\geq \frac{1}{\pi}\lim_{n\to\infty}\|E'_{n+1}\|^2_{L^2(\Omega)}=\frac{1}{2\pi^2}. 
\end{equation}
\end{thm}

The inequality in (\ref{equation2}) is, of course, a consequence of (\ref{equation120}), but it takes some effort to establish the existence and value of the limit in (\ref{equation2}). 

Based on Theorem \ref{thm1}, we find plausible that the following conjecture be true. 
\begin{conj} For an arbitrary piecewise analytic Jordan curve $L$ having at least one corner with exterior angle different from $0$, $\pi$, and $2\pi$, we have
\[  
\lim_{n\to\infty}\|E'_{n+1}\|^2_{L^2(\Omega)}>0.
\]
\end{conj}

The weaker thesis that $\liminf_{n\to\infty}\|E'_{n+1}\|^2_{L^2(\Omega)}>0$ would be enough to guarantee that the $O(1/n)$ estimate for $\alpha_n$ is sharp for every such curve. If the conjecture were true,  it would be interesting to determine whether the value of the limit is, indeed, independent of the curve $L$, and therefore equal to $(2\pi)^{-1}$.

%%%%%%%%%%%%%%%%%%%%%%%%%%%%%%%%%%%%%%%%%%

\section{Proof of Theorem \ref{thm1}}

Hereafter, $L$ will denote the curve in $(\ref{equation1})$.
The other two arcs of the circles $C_1$ and $C_2$ also form a piecewise analytic Jordan curve that we denote by  
\[
\mathcal{L}:=\{z\in C_1: \Im(z)\leq 0\}\cup \{z\in C_2: \Im(z)\geq 0\}.
\]
The exterior of $L$ will be denoted by $\Omega$, while the interior of $\mathcal{L}$ will be denoted by $\mathcal{R}$. 
It is easy to verify that 
\[
\mathcal{R}=\{1/z:z\in \Omega\},
\]
and that the Zhoukowsky transformation $\phi(z)=2^{-1}\left(z+1/z\right)$ maps $\Omega$ conformally onto $\{w:|w|>1\}$. This same function $\phi$ takes both $L$ and $\mathcal{L}$ onto the unit circle.

 With the notation we previously introduced for the Faber polynomials and related quantities, we then have for the curve $L$ that  
\begin{equation}\label{equation4}
\phi^{n}(z)=\frac{1}{2^{n}}\left(z+\frac{1}{z}\right)^{n}=F_{n}(z)+E_{n}(z),
\end{equation}
where $F_{n}$ is the polynomial part of $\phi^{n}$, so that if we define 
\[
G_n(z):=F_n(z)-F_n(0),\quad n\geq 0,
\]
then $G_0(z)\equiv 0$ and
 \begin{equation}\label{explicitGn}
G_{n}(z)=\frac{1}{2^{n}}\sum_{j=0}^{\lfloor \frac{n-1}{2}\rfloor} \binom{n}{j}z^{n-2j}, \quad n\geq 1.
\end{equation}

Since $\phi$ is invariant under $z\mapsto 1/z$, we get from (\ref{equation4}) that 
\begin{align}\label{rel1}
E_{n+1}(z)=G_{n}(1/z),
\end{align}
and that $G_n$ satisfies the recurrence relation 
\begin{align}\label{rel5}
 G_{n+1}(z)={}&\phi(z)G_n(z)+\frac{za_n}{2}-\frac{a_{n+1}}{2},\quad n\geq 0,
\end{align} 
where
\[
a_n:=F_n(0)=\begin{cases}0, & n\ \mathrm{odd},\\
2^{-n}\binom{n}{n/2}, & n\ \mathrm{even},
\end{cases}\quad n\geq 0.
\]
From this explicit expression for $a_n$ one can easily verify that
\begin{align}\label{equation5}
 a_{n}={}& \frac{n-1}{n}a_{n-2},\quad n\geq 2.
\end{align}

At some point, we will need to deal with the quantities
\[
b_{n}:=\int_{-1}^1\frac{G_{n+1}(x)}{x}dx\quad n\geq 0.
\]
If $n$ is odd, $G_{n+1}$ is even and so $b_n=0$. If $n$ is even, then (\ref{explicitGn}) yields
\[
b_{n}=\frac{1}{2^{n}}\sum_{j=0}^{n/2}\frac{\binom{n+1}{j}}{n+1-2j},\quad n=2k,\quad k\geq 0.
\]
From this last expression, it is not difficult to see that 
\begin{align}\label{equation6}
b_{n}={}&\frac{n}{n+1}b_{n-2}+\frac{a_{n}}{n+1}\,,\quad n\geq 2.
\end{align}
Combining (\ref{equation5}) and (\ref{equation6}), we find $(n+2)a_{n+2}b_{n}-na_{n}b_{n-2}=a^2_{n}$, so that 
\begin{equation}\label{equation19}
\sum_{j=0}^ka^2_{2j}=(2k+1)a_{2k}b_{2k}, \quad k\geq 0.
\end{equation}

We now have everything we need to give the proof of Theorem \ref{thm1}.

It follows from (\ref{rel1}), and the fact that $z\mapsto 1/z$ takes $\Omega$ onto the region $\mathcal{R}$, that
\begin{align*}
 \|E'_{n+1}\|^2_{L^2(\Omega)}=& \int_{\mathcal{R}}|G'_{n+1}(z)|^2dxdy.
\end{align*} 

Let $\mathcal{L}_1$ denote the part of $\mathcal{L}$ lying in the closed upper half plane. Since $G_n(\overline{z})=\overline{G_n(z)}$, and since $G'_{n+1}G_{n+1}$ is an odd function, the complex version of Green's formula yields
\begin{align}\label{equation7}
\|E'_{n+1}\|^2_{L^2(\Omega)}={}& I_{n+1,n+1},
\end{align} 
where
\[
I_{n,k}:=\frac{1}{i}\int_{ \mathcal{L}_1}[G_{k}(z)\phi^{n-k}(z)]'\,\overline{G_{k}(z)\phi^{n-k}(z)}dz,\quad n\geq 0,\ 0\leq k\leq n.
\]
Notice that $I_{n,0}=0$. 

Using the recurrence relation (\ref{rel5}), we find that 
\begin{align}\label{rel12}
I_{n+1,k+1}={}&I_{n+1,k}+A_{n,k}+B_{n,k}+C_{n,k},
\end{align} 
where
\begin{align}\label{equation8}
A_{n,k}:={}&\frac{1}{2i}\int_{ \mathcal{L}_1}[\phi^{n+1-k}(z)G_{k}(z)]'\overline{\phi^{n-k}(z)}\left[\overline{z}a_{k}-a_{k+1} \right]dz,
\end{align} 
\begin{align}\label{equation9}
B_{n,k}:={}&\frac{1}{2i}\int_{ \mathcal{L}_1}\overline{\phi^{n+1-k}(z)G_{k}(z)}\left(\phi^{n-k}(z)
\left[za_{k}-a_{k+1}\right]\right)'dz,
\end{align} 
and
\begin{align}\label{equation22}
C_{n,k}:={}&\frac{1}{4i}\int_{ \mathcal{L}_1}\overline{\phi^{n-k}(z)}\left[\overline{z}a_{k}-a_{k+1} \right]\left(\phi^{n-k}(z)\left[za_{k}-a_{k+1}\right]\right)'dz\,.
\end{align} 

We now observe that for $z\in\mathcal{L}_1$, 
\begin{equation}\label{equation10}
|\phi(z)|=1, \quad [\phi(z)]^2\overline{\phi'(z)}\,\overline{dz}=-\phi'(z)dz,
\end{equation}
and
\begin{equation}\label{rel11}
\overline{z}=\frac{1+iz}{z+i}=\frac{2z+i(z^2-1)}{z^2+1},\quad \overline{dz}=-\frac{2dz}{(z+i)^2}=d\left(\frac{1+iz}{z+i}\right).
\end{equation}
Hence, integration by parts in (\ref{equation8}) gives 
\begin{align*}
A_{n,k}={}& \frac{n-k}{2i}\int_{ \mathcal{L}_1}G_{k}(z)\phi'(z)\left[\overline{z}a_{k}-a_{k+1}\right]dz +\frac{a_{k}}{i}\int_{ \mathcal{L}_1}\frac{\phi(z)G_{k}(z)}{(z+i)^2}dz.
\end{align*}

Similarly, by expanding the derivative in (\ref{equation9}), multiplying, and using (\ref{equation10}) and (\ref{rel11}), we obtain that $B_{n,k}=\overline{A_{n,k}}$. Hence,
\begin{align}\label{equation14}
A_{n,k}+B_{n,k}={}& (n-k)\Im\left(\int_{ \mathcal{L}_1}G_{k}(z)\phi'(z)\left[\overline{z}a_{k}-a_{k+1} \right]dz\right)\nonumber\\
{}& + 2a_{k}\Im\left(\int_{ \mathcal{L}_1}\frac{\phi(z)G_{k}(z)}{(z+i)^2}dz\right)\nonumber\\
={}& 2(n+1-k)a_{k}\int_{ -1}^1\frac{G_{k}(z)}{z^2+1}dz- \frac{(n-k)a_{k}}{2}\int_{-1}^1\frac{G_{k}(z)(z^2+1)}{z^2}dz\nonumber\\
{}&+\frac{\pi (n-k)a^2_{k+1}}{2}.
\end{align}
In the calculations leading to the equality of these two last expressions, we have used (\ref{rel11}), that
\[
\frac{(z^2-1)^2}{z^2(z^2+1)} =\frac{z^2+1}{z^2}-\frac{4}{z^2+1},
\]
and that for $k$ even
 \begin{align*}
\int_{-1}^1\frac{G_{k}(z)}{z}\frac{(z^2-1)}{z^2+1}dz=0,
 \end{align*}
while for $k$ odd,
\begin{align*}
\int_{ \mathcal{L}_1}\frac{G_{k}(z)(z^2-1)}{z^2}dz={}&-G'_{k}(0)\int_{ \mathcal{L}_1}\frac{1}{z}dz=-i\pi a_{k+1}.
\end{align*}

Applying the recurrence relation (\ref{rel5}) one more time yields
\begin{align}\label{equation16}
\int_{ -1}^1\frac{G_{k}(z)}{z^2+1}dz={}& \frac{b_{k-2}}{2}-\frac{\pi a_{k}}{4},\qquad \int_{-1}^1\frac{G_{k}(z)(z^2+1)}{z^2}dz=2b_k-2a_{k},
\end{align}
and since $G_0(z)\equiv 0$ and $a_1=0$, we get from (\ref{equation14}), (\ref{equation16}) and (\ref{equation5}) that  
\begin{align}\label{equation18}
\sum_{k=0}^n[A_{n,k}+B_{n,k}]={}& \sum_{k=0}^{n-1}(n-k-1)a_{k+2}b_{k}- \sum_{k=2}^{n}(n-k)a_{k}b_k+\sum_{k=1}^{n}(n-k)a^2_{k}\nonumber\\
={} &n a_{0}b_{0} -(n+1)\sum_{k=0}^{n-1}\frac{a_kb_{k}}{k+2}+\sum_{k=1}^{n}(n-k)a^2_{k}.
\end{align}

We now consider (\ref{equation22}) and compute  
\begin{align*}
C_{n,k}={}&\frac{(n-k)a^2_{k}}{4i}\int_{ \mathcal{L}_1}\overline{\phi(z)}\phi'(z)|z|^2dz+\frac{(n-k)a^2_{k+1}}{4i}\int_{ \mathcal{L}_1}\overline{\phi(z)}\phi'(z)dz\\
&+\frac{a^2_{k}}{4i}\int_{ \mathcal{L}_1}\overline{z}dz\,.
\end{align*}

Using (\ref{equation10}) again, and having in mind that $\phi$ takes $\mathcal{L}_1$ onto the lower half of the unit circle, we get 
 \begin{align*}
\int_{ \mathcal{L}_1}\overline{\phi(z)}\phi'(z)dz={}&\int_{ \mathcal{L}_1}\frac{\phi'(z)}{\phi(z)}dz=-i\pi,
\end{align*}
 \begin{align*}
\int_{ \mathcal{L}_1}\overline{z}dz={}&-\int_{-1}^1\frac{2z}{z^2+1}dz-i\int_{-1}^1\frac{z^2-1}{z^2+1}dz= i(\pi-2)
\end{align*}
and since $|z|^2=1+iz-i\overline{z}$ for $z\in \mathcal{L}_1$,
\begin{align*}
\int_{ \mathcal{L}_1}\overline{\phi(z)}\phi'(z)|z|^2dz={}&-i\pi+ i\int_{ \mathcal{L}_1}\frac{\phi'(z)}{\phi(z)}zdz+ i\overline{\int_{\mathcal{L}_1}\frac{\phi'(z)}{\phi(z)}zdz}\\
={}&-i\pi-2i\int_{-1}^{1}\frac{z^2-1}{z^2+1}dz.
\end{align*}
Hence,
\begin{align}\label{equation17}
C_{n,k}={}&\frac{\pi}{4}\left[(n+1-k)a^2_{k}-(n-k)a^2_{k+1}\right]-(n-k)a^2_k-a^2_{k}/2.
\end{align}

Then, combining (\ref{rel12}), (\ref{equation18}) and (\ref{equation17}), and since $a_0=b_0=1$, we obtain 
\begin{align*}
I_{n+1,n+1}={}&\sum_{k=0}^n\left(A_{n,k}+B_{n,k}+C_{n,k}\right)=  (n+1)\left[\frac{\pi}{4}-\sum_{k=0}^{n-1}\frac{a_kb_{k}}{k+2}\right]-\frac{1}{2}\sum_{k=0}^{n}a^2_{k}.
\end{align*}

Now, summation by parts gives
\begin{align*}
\frac{1}{2N+1}\sum_{k=0}^Na^2_{2k}={}&\sum_{k=0}^N\frac{a^2_{2k}}{2k+1}-\sum_{k=0}^{N-1}\left(\frac{1}{2k+1}-\frac{1}{2k+3}\right)\sum_{j=0}^ka^2_{2j}.
\end{align*}

Combining the two last equalities and (\ref{equation19}), we get that for every integer $N\geq 0$, we have 
\begin{align}\label{equation20}
\frac{I_{2N+1,2N+1}}{2N+1}={}&\frac{\pi}{4}-\frac{1}{2}\sum_{k=0}^N\frac{a^2_{2k}}{2k+1}-\sum_{k=0}^{N-1}\frac{\sum_{j=0}^ka^2_{2j}}{(2k+1)(2k+2)(2k+3)},
\end{align}
\begin{align*}
I_{2N+2,2N+2}={}& \frac{(2N+2)}{(2N+1)}I_{2N+1,2N+1}-\frac{\sum_{j=0}^Na^2_{2j}}{(2N+1)}.
\end{align*}

Since (\ref{equation121}) tells us that $I_{n+1,n+1}$, which is defined by (\ref{equation7}), is bounded above, the bracket in (\ref{equation20}) must converge to zero as $n\to\infty$, and since 
\begin{align*}
a_{2k}=\frac{1}{2^{2k}}&\binom{2k}{k}=\frac{\Gamma(1/2)\Gamma(k+1/2)}{\pi\Gamma(k+1)}=\frac{1+O(1/k)}{\sqrt{\pi (k-1/2)}},
\end{align*}
we arrive at   
\begin{align*}
\frac{I_{2N+1,2N+1}}{2N+1}={}& \frac{1}{2}\sum_{k=N+1}^\infty\frac{a^2_{2k}}{2k+1}+\sum_{k=N}^{\infty}\frac{\sum_{j=0}^ka^2_{2j}}{(2k+1)(2k+2)(2k+3)}\\
 ={} & \frac{1}{2\pi}\sum_{k=N+1}^\infty\frac{2+O(1/k)}{(2k-1)(2k+1)}+\sum_{k=N}^{\infty}\frac{O(\sum_{j=0}^k\frac{1}{j})}{(2k+1)(2k+2)(2k+3)}\\
={}& \frac{1}{2\pi(2N+1)}+O(\ln N/N^2),
\end{align*}
and Theorem \ref{thm1} follows.


\begin{thebibliography}{99}
\bibitem{Carleman}{T. Carleman, \"{U}ber die approximation analytischer funktionen durch lineare aggregate von vorgegebenen
potenzen,  Archiv. f\"{o}r Math. Atron. och Fysik, 17
(1922) 1-30.}

\bibitem{DMN}{P. Dragnev, E. Mi\~{n}a-D\'{\i}az, Michael Northington V, Asymptotics of Carleman polynomials for level curves of a shifted Zhukovsky transformation, Comput. Methods Funct. Theory 13 (2013) 75-89. }

\bibitem{DM1}{P. Dragnev, E. Mi\~{n}a-D\'{\i}az, On a series representation for Carleman orthogonal polynomials, Proc. Amer. Math. Soc. 138 (2010) 4271-4279.}

\bibitem{DM2}{P. Dragnev, E. Mi\~{n}a-D\'{\i}az,  Asymptotic behavior and zero distribution of Carleman orthogonal polynomials,  J. Approx. Theory. 162 (2010) 1982-2003}

\bibitem{Levin}{A. L. Levin, E. B. Saff, N. S. Stylianopoulos,
Zero distribution of Bergman orthogonal polynomials for certain
planar domains, Constr. Approx., 19 (2003), 411-435.}

\bibitem{May}{V. Maymeskul, E. B. Saff, Zeros of polynomials orthogonal over
regular $N$-gons, J. Approx. Theory, 122 (2003), 129-140.}

\bibitem{M1}{E. Mi\~{n}a-D\'{\i}az, An asymptotic integral representation
for Carleman orthogonal polynomials, Int. Math. Res. Notices,  2008 (2008), article ID rnn065, 38 pages.}

\bibitem{MSS}{E. Mi\~{n}a-D\'{\i}az,  E. B. Saff, N. S. Stylianopoulos, Zero distributions
for polynomials orthogonal with weights over
certain planar regions, Comput. Methods Funct. Theory, 5 (2005), 185-221.}

\bibitem{NS}{N. Stylianopoulos, Strong asymptotics of Bergman polynomials over domains with corners and applications, Constr. Approx. 38 (2013), 59-100.} 

\bibitem{Suetin}{P. K. Suetin, Polynomials orthogonal over a region and Bieberbach polynomials, vol. 100 of Proc. Steklov Inst. Math., Amer. Math. Soc. Translations, 1974.}

\bibitem{Suetin2}{P. K. Suetin, Series of Faber Polynomials, Gordon and Breach Science Publications, Amsterdam, 1998.}
\end{thebibliography}
\end{document}